\documentclass[11pt]{amsart}
\usepackage{amscd,amssymb}
\usepackage{amsthm,amsmath,amssymb}
\usepackage[matrix,arrow]{xy}

\newtheorem{theorem}[equation]{Theorem}
\newtheorem{proposition}[equation]{Proposition}
\newtheorem{lemma}[equation]{Lemma}
\newtheorem{corollary}[equation]{Corollary}

\theoremstyle{definition}

\theoremstyle{remark}

\makeatletter\@addtoreset{equation}{section} \makeatother

\author{Ivan Cheltsov}

\title{Nonrational del Pezzo fibrations}

\thanks{The author would like to thank A.\,Corti, M.\,Grinenko,
V.\,Is\-kov\-skikh, V.\,Sho\-ku\-rov for fruitful conversations.}

\begin{document}

\begin{abstract}
Let $X$~be a general~divisor~in $|3M+nL|$~on the~rational scroll
$\mathrm{Proj}(\oplus_{i=1}^{4}\mathcal{O}_{\mathbb{P}^{1}}(d_{i}))$,
where $d_{i}$ and $n$ are~integers, $M$ is the~tautological line
bundle, $L$ is a fibre of the~natural projection to
$\mathbb{P}^{1}$, and $d_{1}\geqslant\cdots\geqslant d_{4}=0$. We
prove that $X$ is rational $\iff$ $d_{1}=0$~and~$n=1$.
\end{abstract}

\maketitle

\section{Introduction.}
\label{section:introduction}

The rationality problem for threefolds\footnote{All varieties are
assumed to be projective, normal, and defined over $\mathbb{C}$.}
splits in three cases: conic bundles, del Pezzo fibrations, and
Fano threefolds. The~cases of conic bundles and Fano threefolds
are well studied.

Let $\psi\colon X\to\mathbb{P}^{1}$ be a fibration into del Pezzo
surfaces of degree $k\geqslant 1$ such that $X$ is smooth and
$\mathrm{rk}\,\mathrm{Pic}(X)=2$. Then $X$ is rational if
$k\geqslant 5$.  The~following result is due to \cite{Al87} and
\cite{Shr05}.

\begin{theorem}
\label{theorem:Valera} Suppose that fibres of $\psi$ are normal
and $k=4$. Then $X$ is rational if~and~only~if
$$
\chi\big(X\big)\in\big\{0, -8, -4\big\},
$$
where $\chi(X)$ is the~topological
Euler characteristic.
\end{theorem}

The following result is due to \cite{Pu98}.

\begin{theorem}
\label{theorem:Sasha-I} Suppose that $K_{X}^{2}\not\in
\mathrm{Int}\,\overline{\mathrm{NE}}(X)$ and $k\leqslant 2$. Then
$X$ is nonrational.
\end{theorem}

In the~case when $k\leqslant 2$ and $K_{X}^{2}\in
\mathrm{Int}\,\overline{\mathrm{NE}}(X)$, the~threefold $X$
belongs to finitely many~deformation families, whose general
members are nonrational (see \cite{Vo88}, \cite{Ch04b},
\cite{Gr05}, Proposition~\ref{proposition:last-dp2}).

Suppose that $k=3$. Then $X$ is a divisor in the~linear system
$|3M+nL|$ on the~scroll
$$
\mathrm{Proj}\Big(\oplus_{i=1}^{4}\mathcal{O}_{\mathbb{P}^{1}}\big(d_{i}\big)\Big),
$$
where $n$ and $d_{i}$ are integers, $M$ is the~tautological line
bundle, and $L$ is a fibre of the~natural projection to
$\mathbb{P}^{1}$. Suppose that  $d_{1}\geqslant d_{2}\geqslant
d_{3}\geqslant d_{4}=0$.

Suppose that $X$ is a general\footnote{A complement to a countable
union of Zariski closed subsets.} divisor in $|3M+nL|$. The
following result is due to \cite{Pu98}.

\begin{theorem}
\label{theorem:Sasha-II} Suppose that $K_{X}^{2}\not\in
\mathrm{Int}\,\overline{\mathrm{NE}}(X)$. Then $X$ is nonrational.
\end{theorem}

It follows from \cite{ClGr72}, \cite{Sho83}, \cite{Bar84},
\cite{Vo88}, \cite{BrCoZu03}, \cite{Ch04b} that $X$ is nonrational
when
$$
\big(d_{1},d_{2},d_{3},n\big)\in\Big\{\big(0,0,0,2\big),\big(1,0,0,0\big),\big(2,1,1,-2\big),\big(1,1,1,-1\big)\Big\}.%
$$

We prove the~following~result in Section~\ref{section:main}.

\begin{theorem}
\label{theorem:main} The threefold $X$ is rational $\iff$
$d_{1}=0$ and $n=1$.
\end{theorem}

Therefore, the~threefold $X$ is nonrational in the case when
$\chi(X)\ne -14$. Indeed, we have
$$
\chi(X)=-4K_{X}^{3}-54=-4\Big(18-6\big(d_{1}+d_{2}+d_{3}\big)-8n\Big)-54=18-24\big(d_{1}+d_{2}+d_{3}\big)-32n,
$$
and $\chi(X)=-14$ implies that $(d_{1}, d_{2}, d_{3}, n)=(0, 0, 0,
1)$ or $(d_{1}, d_{2}, d_{3}, n)=(2, 1, 1, -2)$.

The inequality $5n\geqslant 12-3(d_{1}+d_{2}+d_{3})$ holds when
$K_{X}^{2}\not\in \mathrm{Int}\,\overline{\mathrm{NE}}(X)$. For
$n<0$,~the~inequality
$$
5n\geqslant 12-3\big(d_{1}+d_{2}+d_{3}\big)
$$
implies that $K_{X}^{2}\not\in
\mathrm{Int}\,\overline{\mathrm{NE}}(X)$ (see Lemma~36 in
\cite{BrCoZu03}). Hence, the threefold $X$ does not~belong~to
finitely many deformation families in the~case when
$K_{X}^{2}\in\mathrm{Int}\,\overline{\mathrm{NE}}(X)$ (see
Section~\ref{section:preliminaries}).

Let us illustrate our methods by proving the~following result.

\begin{proposition}
\label{proposition:last-dp2} Let $X$ be double cover of the scroll
$$
\mathrm{Proj}\Big(\mathcal{O}_{\mathbb{P}^{1}}\big(2\big)\oplus\mathcal{O}_{\mathbb{P}^{1}}(2)\oplus\mathcal{O}_{\mathbb{P}^{1}}\Big)
$$
that is branched over a general\footnote{A complement to a
countable union of Zariski closed subsets.} divisor $D\in
|4M-2L|$, where $M$ is the~tautological~line bundle,
and~$L$~is~a~fibre of the~natural projection to $\mathbb{P}^{1}$.
Then $X$ is nonrational.
\end{proposition}

\begin{proof}
Put
$V=\mathrm{Proj}(\mathcal{O}_{\mathbb{P}^{1}}(2)\oplus\mathcal{O}_{\mathbb{P}^{1}}(2)\oplus\mathcal{O}_{\mathbb{P}^{1}})$.
The~divisor $D$ is given by the equation
\begin{multline*}
\alpha_6x_1^4+\alpha^1_6x_1^3x_2+\alpha_4x_1^3x_3+\alpha^2_6x_1^2x_2^2+\alpha^1_4x_1^2x_2x_3+\alpha_2x_1^2x_3^2+\alpha^3_6x_1x_2^3+\\
+\alpha^2_4x_1x_2^2x_3+\alpha^1_2x_1x_2x_3^2+\alpha_0x_1x_3^3+\alpha^4_6x_2^4+\alpha^3_4x_2^3x_3+\alpha^2_2x_2^2x_3^2+\alpha^1_0x_2x_3^3=0
\end{multline*}
in bihomogeneous coordinates on $V$ (see \S 2.2 in \cite{Re97}),
where $\alpha_d^{i}=\alpha_d^{i}(t_1,t_2)$ is a sufficiently
general homogeneous polynomial of degree $d\geqslant 0$. Let
$$
\chi\colon Y\longrightarrow
\mathrm{Proj}\Big(\mathcal{O}_{\mathbb{P}^{1}}\big(2\big)\oplus\mathcal{O}_{\mathbb{P}^{1}}(2)\oplus\mathcal{O}_{\mathbb{P}^{1}}\Big)
$$
be a double cover branched over a divisor $\Delta\subset V$ that
is given by the~same bihomogeneous equation as of divisor $D$ with
the~only exception that $\alpha_0=\alpha_0^{1}=0$. Then $Y$ is
singular, because the~divisor $\Delta$ is singular along the~curve
$Y_{3}\subset V$ that is given by the~equations $x_1=x_2=0$.

The Bertini theorem implies the~smoothness of $\Delta$ outside of
the~curve $Y_{3}$.

Let $C$ be a curve on the~threefold $Y$ such that $\chi(C)=Y_{3}$.
Then the threefold $Y$ has~sin\-gu\-la\-ri\-ties of type
$\mathbb{A}_{1}\times\mathbb{C}$ at general point of the~curve
$C$. We may assume that the system
$$
\alpha_2\big(t_1,t_2\big)=\alpha^1_2\big(t_1,t_2\big)=\alpha^2_2\big(t_1,t_2\big)=0
$$
has no non-trivial solutions. Then $Y$ has singularities of type
$\mathbb{A}_{1}\times\mathbb{C}$ at every point of $C$.

Let $\alpha\colon\tilde{V}\to V$ be the~blow up of $Y_{3}$, and
$\beta\colon\tilde{Y}\to Y$ be the~blow up of $C$. Then
the~diagram
$$
\xymatrix{
\tilde{Y}\ar@{->}[d]_{\beta}\ar@{->}[rr]^{\tilde{\chi}}&& \tilde{V}\ar[d]^{\alpha}\\
Y\ar@{->}[rr]^{\chi}&&V}
$$
commutes, where $\tilde{\chi}\colon\tilde{Y}\to\tilde{V}$ is a
double cover. The threefold $\tilde{Y}$ is smooth.

Let $E$ be the~exceptional divisor of $\alpha$, and
$\tilde{\Delta}$ be the proper transform of $\Delta$ via $\alpha$.
Then
$$
\tilde{\Delta}\sim \alpha^{*}\big(4M-2L\big)-2E,%
$$
which implies that $\tilde{\Delta}$ is nef and big, because the
pencil $|\alpha^{*}(M-2L)-E|$ does not have base points. The
morphism $\tilde{\chi}$ is branched over $\tilde{\Delta}$. Then
$\mathrm{rk}\,\mathrm{Pic}(\tilde{Y})=3$ by Theorem~2 in
\cite{Ra05}.

The linear system $|g^{*}(M-L)-E|$ does not have base points and
gives a $\mathbb{P}^{1}$-bundle
$$
\tau\colon\tilde{V}\longrightarrow\mathrm{Proj}\Big(\mathcal{O}_{\mathbb{P}^{1}}\big(2\big)\oplus\mathcal{O}_{\mathbb{P}^{1}}\big(2\big)\Big)\cong\mathbb{F}_{0},
$$
which induces a conic bundle
$\tilde{\tau}=\tau\circ\tilde{\chi}\colon
\tilde{Y}\to\mathbb{F}_{0}$.

Let $Y_{2}\subset V$ be the~subscroll given by $x_{1}=0$, and $S$
be a proper transform of $Y_{2}$ via $\alpha$. Then
$$
Y_{2}\cong\mathrm{Proj}\Big(\mathcal{O}_{\mathbb{P}^{1}}\big(2\big)\oplus\mathcal{O}_{\mathbb{P}^{1}}\Big)\cong\mathbb{F}_{2},
$$
and $S\cong Y_{2}$. But $\tau$ maps $S$ to the~section of
$\mathbb{F}_{0}$ that has trivial self-intersection.

Let $\tilde{S}$ be a surface in $\tilde{Y}$ such that
$\tilde{\chi}(\tilde{S})=S$, and $Z\subset\tilde{Y}$ be a general
fibre of the~natural projection to $\mathbb{P}^{1}$. Then $-K_{Z}$
is nef and big and $K_{Z}^{2}=2$. But the~morphism
$$
\alpha\circ\tilde{\chi}\big\vert_{\tilde{S}}\colon\tilde{S}\longrightarrow Y_{2}%
$$
is a double cover branched over a divisor that is cut out by the
equation
$$
\alpha^4_6\big(t_{0},t_{1}\big)x_2^2+\alpha^3_4\big(t_{0},t_{1}\big)x_2x_3+\alpha^2_2\big(t_{0},t_{1}\big)x_3^2=0.
$$

Let $\Xi\subset\mathbb{F}_{0}$ be a degeneration divisor of the
conic bundle $\tilde{\tau}$. Then
$$
\Xi\sim \lambda\tilde{\tau}\big(\tilde{S}\big)+\mu\tilde{\tau}\big(Z\big),%
$$
where $\lambda$ and $\mu$ are integers. But $\lambda=6$, because
$K_{Z}^{2}=2$. We have $\tilde{\tau}(\tilde{S})\not\subset\Xi$.
Then
$$
\mu=\tilde{\tau}(\tilde{S})\cdot\Xi=8-K_{\tilde{S}}^{2},
$$
because $\mu$ is the number of reducible fibres of the conic
bundle $\tilde{\tau}\vert_{\tilde{S}}$. These fibers are given by
$$
\Big(\alpha^3_4\big(t_{0},t_{1}\big)\Big)^{2}=4\alpha^2_2\big(t_{0},t_{1}\big)\alpha^4_6\big(t_{0},t_{1}\big),
$$
which implies that $\mu=\tilde{\tau}(\tilde{S})\cdot \Xi=8$. Then
$\tilde{Y}$ is nonruled by Theorem~10.2 in \cite{Sho83}, which
implies the nonrationality of the threefold $X$ by Theorem~1.8.3
in \S IV of the~book \cite{Ko96}.
\end{proof}

\section{Preliminaries.}
\label{section:preliminaries}

All results of this section follow from \cite{BrCoZu03}. Take a
scroll
$$
V=\mathrm{Proj}\Big(\oplus_{i=1}^{4}\mathcal{O}_{\mathbb{P}^{1}}\big(d_{i}\big)\Big),
$$
where $d_{i}$ is an integer, and $d_{1}\geqslant d_{2}\geqslant
d_{3}\geqslant d_{4}=0$. Let $M$ and $L$ be the~tautological line
bundle and a fibre of the~natural projection to $\mathbb{P}^{1}$,
respectively. Then $\mathrm{Pic}(V)=\mathbb{Z}M\oplus\mathbb{Z}L$.

Let $(t_{1}:t_{2};x_1:x_{2}:x_{3}:x_k)$ be bihomogeneous
coordinates on $V$ such that  $x_{i}=0$ defines a~divisor in
$|M-d_{i}L|$, and $L$ is given by $t_{1}=0$. Then $|aM+bL|$ is
spanned by divisors
$$
c_{i_{1}i_{2}i_{3}i_{4}}\big(t_{1},t_{2}\big)x_1^{i_{1}}
x_2^{i_{2}}x_3^{i_{3}}x_k^{i_{4}}=0,
$$
where $\sum_{j=1}^{4}i_{j}=a$ and
$c_{i_{1}i_{2}i_{3}i_{4}}(t_{1},t_{2})$ is a homogeneous
polynomial of degree $b+\sum_{j=1}^{4}i_{j}d_{j}$.

Let $Y_{j}\subseteq V$ be a subscroll $x_{1}=\cdots=x_{j-1}=0$.
The~following result holds (see \S 2.8 in \cite{Re97}).

\begin{corollary}
\label{corollary:lemma-of-Reid} Take $D\in |aM+bL|$ and
$q\in\mathbb{N}$, where $a$ and $b$ are integers. Then
$$
\mathrm{mult}_{Y_{j}}\big(D\big)\geqslant q\iff
ad_{j}+b+\big(d_{1}-d_{j}\big)\big(q-1\big)<0.
$$
\end{corollary}

Let $X$ be a general\footnote{A complement to a Zariski closed
subset in moduli.} divisor in $|3M+nL|$, where $n$ is an integer.

\begin{lemma}
\label{lemma:smoothness} Suppose $X$ is smooth and
$\mathrm{rk}\,\mathrm{Pic}(X)=2$. Then $d_{1}\geqslant -n$ and
$3d_{3}\geqslant -n$.
\end{lemma}

\begin{proof}
We see that $Y_{2}\not\subset X$. Then $Y_{3}\not\subset X$,
because $\mathrm{rk}\,\mathrm{Pic}(X)=2$. But
$\mathrm{mult}_{Y_{4}}\big(X\big)\leqslant 1$, because the
threefold $X$ is smooth. The assertion of
Corollary~\ref{corollary:lemma-of-Reid} concludes the~proof.
\end{proof}

\begin{lemma}
\label{lemma:smoothness-extra} Suppose $X$ is smooth and
$\mathrm{rk}\,\mathrm{Pic}(X)=2$. Then either $d_{1}=-n$ or
$d_{2}\geqslant -n$.
\end{lemma}

\begin{proof}
Suppose that $r=d_{1}+n>0$ and $d_{2}<-n$. Then $X$ can be given
by the~equation
$$
\sum_{i,\,j,\,k\geqslant 0\atop
i+j+k=2}\gamma_{ijk}(t_{0},t_{2})x_1^{i}x_{2}^{j}x_{3}^{k}x_4=\alpha_{r}(t_1,t_2)x_{1}x_{4}^{2}+\sum_{i,\,j,\,k\geqslant 0\atop i+j+k=3}\beta_{ijk}(t_{0},t_{2})x_1^{i}x_{2}^{j}x_{3}^{k},%
$$
where $\alpha_{r}(t_1,t_2)$ is a homogeneous polynomial of degree
$r$, $\beta_{ijk}$ and $\gamma_{ijk}$ are homogeneous polynomial
of degree $n+id_{1}+jd_{2}+kd_{3}$. Then every point of
the~intersection
$$
x_{1}=x_{2}=x_{3}=\alpha_{r}\big(t_1,t_2\big)=0
$$
must be singular on the~threefold $X$, which is a contradiction.
\end{proof}

\begin{lemma}
\label{lemma:smoothness-rkPic-2} Suppose $X$ is smooth,
$d_{2}=d_{3}$, $n<0$ and $\mathrm{rk}\,\mathrm{Pic}(X)=2$. Then
$3d_{3}\ne -n$.
\end{lemma}

\begin{proof}
Suppose that $3d_{3}=-n$. Then $X$ can be given by
the~the~bihomogeneous equation
$$
\sum_{j,\,k,\,l\geqslant 0\atop i+j+k=2}\gamma_{jkl}(t_{0},t_{2})x_1x_{2}^{j}x_{3}^{k}x_{4}^{l}=f_{3}(x_{2},x_{3})+\alpha_{r}(t_{0},t_{2})x_{1}^{3}+\sum_{j,\,k,\,l\geqslant 0\atop j+k+l=1}\beta_{jkl}(t_{0},t_{2})x_1^{2}x_{2}^{j}x_{3}^{k}x_{4}^{l},%
$$
where $f_{3}(x_{2},x_{3})$ is a homogeneous polynomial of degree
$3$, $\beta_{jkl}$ and $\gamma_{jkl}$ are homogeneous polynomial
of degree $n+2d_{1}+jd_{2}+kd_{3}$ and $n+d_{1}+jd_{2}+kd_{3}$
respectively, $\alpha_{r}$ is a homogeneous polynomial of degree
$r=3d_{1}+n$. The threefold $X$ contains $3$ subscrolls given by
the~equations
$$
x_{1}=f_{3}(x_{2},x_{3})=0,
$$
which is impossible, because $\mathrm{rk}\,\mathrm{Pic}(X)=2$.
\end{proof}

The following result follows from Lemmas~\ref{lemma:smoothness},
\ref{lemma:smoothness-extra} and \ref{lemma:smoothness-rkPic-2}.

\begin{lemma}
\label{lemma:smoothness-reverse} The threefold $X$ is smooth and
$\mathrm{rk}\,\mathrm{Pic}(X)=2$ whenever
\begin{enumerate}
\item in the~case when $d_{1}=0$, the~inequality $n>0$ holds,%
\item either $d_{1}=-n$ and $3d_{3}\geqslant -n$, or $d_{1}>-n$, $d_{2}\geqslant -n$ and $3d_{3}\geqslant -n$,%
\item in the~case when $d_{2}=d_{3}$ and $n<0$, the~inequality $3d_{3}>-n$ holds.%
\end{enumerate}
\end{lemma}

\begin{proof}
Suppose that all these conditions are satisfied. We must show that
$X$ is smooth, because the~equality
$\mathrm{rk}\,\mathrm{Pic}(X)=2$ holds by Proposition~32 in
\cite{BrCoZu03}.

The linear system $|3M+nL|$ does not have base points if
$n\geqslant 0$. So, the threefold~$X$~is~smooth by the~Bertini
theorem in the~case $n\geqslant 0$. Therefore, we may assume that
$n<0$.

The base locus of $|3M+nL|$ consists of the~curve $Y_{4}$, which
implies that $X$ is smooth outside of the~curve $Y_{4}$ and in a
general point of $Y_{4}$ by the~Bertini theorem and
Corollary~\ref{corollary:lemma-of-Reid}, respectively.

In the~case when $d_{1}=-n$ and $d_{2}<-n$, the~bihomogeneous
equation of the threefold $X$ is
$$
\sum_{i,\,j,\,k\geqslant 0\atop i+j+k=2}\gamma_{ijk}(t_{0},t_{2})x_1^{i}x_{2}^{j}x_{3}^{k}x_4=\alpha_{0}x_{1}x_{4}^{2}+\sum_{i,\,j,\,k\geqslant 0\atop i+j+k=3}\beta_{ijk}(t_{0},t_{2})x_1^{i}x_{2}^{j}x_{3}^{k},%
$$
where $\beta_{ijk}$ and $\gamma_{ijk}$ are homogeneous polynomials
of degree $n+id_{1}+jd_{2}+kd_{3}$ and $\alpha_{0}$ is~a~nonzero
constant. The curve $Y_{4}$ is given by  $x_{1}=x_{2}=x_{3}=0$,
which implies that $X$ is smooth.

In the~case when $d_{1}>-n$ and $d_{2}\geqslant -n$,
the~bihomogeneous equation of $X$ is
$$
\sum_{i,\,j,\,k\geqslant 0\atop i+j+k=2}\gamma_{ijk}(t_{0},t_{2})x_1^{i}x_{2}^{j}x_{3}^{k}x_4=\sum_{i=1}^{3}\alpha_{i}(t_{0},t_{2})x_{i}x_{4}^{2}+\sum_{i,\,j,\,k\geqslant 0\atop i+j+k=3}\beta_{ijk}(t_{0},t_{2})x_1^{i}x_{2}^{j}x_{3}^{k},%
$$
where $\alpha_{i}$ is a homogeneous polynomial of degree
$d_{i}+n$, and $\beta_{ijk}$ and $\gamma_{ijk}$ are homogeneous
polynomials of degree $n+id_{1}+jd_{2}+kd_{3}$. Therefore, tither
$\alpha_{1}x_{1}x_{4}^{2}$ or $\alpha_{2}x_{2}x_{4}^{2}$ does not
vanish at any given point of the~curve $Y_{4}$, which implies that
$X$ is smooth.
\end{proof}

Thus, there is an infinite series of quadruples
$(d_{1},d_{2},d_{3},n)$ such that the threefold $X$ is smooth, the
equality $\mathrm{rk}\,\mathrm{Pic}(X)=2$ holds, the inequality
$5n<12-3(d_{1}+d_{2}+d_{3})$ holds and $n<0$.

\section{Nonrationality.}
\label{section:main}

We use the~notation of Section~\ref{section:preliminaries}. Let
$X$ be a general\footnote{A complement to a countable union of
Zariski closed subsets.} divisor in $|3M+nL|$, and suppose~that
the~threefold $X$ is smooth, $\mathrm{rk}\,\mathrm{Pic}(X)=2$, and
$X$ is rational. Let us show that $d_{1}=0$ and $n=1$.

The threefold $X$ is given by a~bihomogeneous equation
$$
\sum_{l=0}^{3}\alpha_{i}\big(t_{0},t_{2}\big)x_3^{i}x_{4}^{3-i}+x_{1}F\big(t_{0},t_{1},x_{1},x_{2},x_{3},x_{4}\big)+x_{2}G\big(t_{0},t_{1},x_{1},x_{2},x_{3},x_{4}\big)=0,
$$
where $\alpha_{i}$ is a general homogeneous polynomial of degree
$n+id_{3}$, and $F$ and $G$ stand for
$$
\sum_{i,\,j,\,k,\,l\geqslant 0\atop i+j+k+l=2}\beta_{ijkl}(t_{0},t_{2})x_1^{i}x_{2}^{j}x_{3}^{k}x_{4}^{l}\ \ \text{and}\sum_{i,\,j,\,k,\,l\geqslant 0\atop i+j+k+l=2}\gamma_{ijkl}(t_{0},t_{2})x_1^{i}x_{2}^{j}x_{3}^{k}x_{4}^{l}%
$$
respectively, where $\beta_{ijkl}$ is a general homogeneous
polynomial of degree $n+(i+1)d_{1}+jd_{2}+kd_{3}$,
and~$\gamma_{ijkl}$ is a general homogeneous polynomial of degree
$n+id_{1}+(j+1)d_{2}+kd_{3}$.

Let $Y$ be a threefold given by $x_{1}F+x_{2}G=0$. Then
$Y_{3}\subset Y$, where $Y_{3}$ is given by $x_{1}=x_{2}=0$.

\begin{lemma}
\label{lemma:singularities-of-degeneration} The threefold $Y$ has
$2d_{1}+2d_{2}+4d_{3}+4n>0$ isolated ordinary double points.
\end{lemma}

\begin{proof}
The threefold $Y$ is singular exactly at the~points of $V$ where
$$
x_{1}=x_{2}=F\big(t_{0},t_{1},x_{1},x_{2},x_{3},x_{4}\big)=G\big(t_{0},t_{1},x_{1},x_{2},x_{3},x_{4}\big)=0
$$
by the~Bertini theorem. But
$Y_{3}\cong\mathrm{Proj}(\mathcal{O}_{\mathbb{P}^{1}}(d_{3})\oplus\mathcal{O}_{\mathbb{P}^{1}})\cong\mathbb{F}_{d_{3}}$,
where $(t_{0}:t_{1};x_{3}:x_{4})$ can be considered as natural
bihomogeneous coordinates on the~surface $Y_{3}$.

Let $C$ and $Z$ be the~curves on  $Y_{3}$ that are cut out by the
equations $F=0$ and $G=0$,~respectively. Then $C$ and $Z$ are
given by the equations
$$
\sum_{k,\,l\geqslant 0\atop k+l=2}\beta_{kl}(t_{0},t_{2})x_{3}^{k}x_{4}^{l}=0\ \ \text{and}\ \sum_{k,\,l\geqslant 0\atop k+l=2}\gamma_{kl}(t_{0},t_{2})x_{3}^{k}x_{4}^{l}=0%
$$
respectively, where $\beta_{kl}=\beta_{00kl}$ and
$\gamma_{kl}=\gamma_{00kl}$.

The degrees of $\beta_{kl}$ and $\gamma_{kl}$ are $n+d_{1}+kd_{3}$
and $n+d_{2}+kd_{3}$, respectively.

Let $O$ be a point of the~scroll $V$ such that the set
$$
x_{1}=x_{2}=F\big(t_{0},t_{1},x_{1},x_{2},x_{3},x_{4}\big)=G\big(t_{0},t_{1},x_{1},x_{2},x_{3},x_{4}\big)=0
$$
contains the point $O$. Then $O\in C\cap Z$ and
$O\in\mathrm{Sing}(Y)$.

It is easy to see that $O$ is an isolated ordinary double point of
the threefold $Y$ in the case~when the~curves $C$ and $Z$ are
smooth and intersect each other transversally at the~point $O$.

Put $\bar{M}=M\vert_{Y_{3}}$ and $\bar{L}=L\vert_{Y_{3}}$. Then
$C\in |2\bar{M}+(n+d_{1})\bar{L}|$ and
$Z\in|2\bar{M}+(n+d_{2})\bar{L}|$. But
$$
\Big|2\bar{M}+\big(n+d_{1}\big)\bar{L}\Big|
$$
does not have base points, because $d_{1}+n\geqslant 0$ by
Lemma~\ref{lemma:smoothness}. So, the~curve $C$ is smooth.

The linear system $|2\bar{M}+(n+d_{2})\bar{L}|$ may have base
components, and $Z$ may not be reduced or irreducible. We have to
show that $C$ intersects $Z$ transversally at smooth points of
$Z$, because
$$
\big|C\cap Z\big|=C\cdot Z=2d_{1}+2d_{2}+4d_{3}+4n,
$$
where $2d_{1}+2d_{2}+4d_{3}+4n>0$ by
Lemmas~\ref{lemma:smoothness}, \ref{lemma:smoothness-extra} and
\ref{lemma:smoothness-rkPic-2}.

Suppose that $d_{1}>-n$. Then $d_{2}\geqslant -n$ by
Lemma~\ref{lemma:smoothness-extra}. We see that
$|2\bar{M}+(n+d_{2})\bar{L}|$ does not have base points. Then $Z$
is smooth and $C$ intersects $Z$ transversally at every point of
$C\cap Z$.

We may assume that $d_{1}=-n$. Let $Y_{4}\subset Y_{3}$ be a curve
given by $x_{3}=0$. Then
$$
C\cap Y_{4}=\emptyset,
$$
and either the~linear system $|2\bar{M}+(n+d_{2})\bar{L}|$ does
not have base points, or the~base locus of the~linear system
$|2\bar{M}+(n+d_{2})\bar{L}|$ consist of the~curve $Y_{4}$.
However, we have
$$
C\cap Z\subset Y_{3}\setminus Y_{4},
$$
which implies that $C$ intersects the curve $Z$ transversally at
smooth points of $Z$.
\end{proof}

Let $\pi\colon\tilde{V}\to V$ be the~blow up of $Y_{3}$, and
$\tilde{Y}$ be and the~proper transforms of $Y$ via $\pi$. Then
$$
\tilde{Y}\sim \pi^{*}\big(3M+nL\big)-E,
$$
where $E$ is and exceptional divisor of $\pi$. The threefold
$\tilde{Y}$ is smooth.

\begin{lemma}
\label{lemma:class-group-of-degeneration} The equality
$\mathrm{rk}\,\mathrm{Pic}(\tilde{Y})=3$ holds.
\end{lemma}

\begin{proof}
The~linear system $|\pi^{*}(M-d_{2}L)-E|$ does not have base
points. Thus, the~divisor
$$
\tilde{Y}\sim \pi^{*}\big(3M+nL\big)-E
$$
is nef and big when $n\geqslant 0$ by
Lemmas~\ref{lemma:smoothness}, \ref{lemma:smoothness-extra} and
\ref{lemma:smoothness-rkPic-2}. Hence, the equality
$\mathrm{rk}\,\mathrm{Pic}(\tilde{Y})=3$ holds in the~case when
$n\geqslant 0$ by Theorem~2 in \cite{Ra05}. So, we may assume that
$n<0$.

Let $\omega\colon\tilde{Y}\to\mathbb{P}^{1}$ be the~natural
projection and $S$ be the~generic fibre of $\omega$, which is
considered as a surface defined over the~function field
$\mathbb{C}(t)$. Then $S$ is a smooth cubic surface in
$\mathbb{P}^{3}$, which contains a line in $\mathbb{P}^{3}$
defined over the~field $\mathbb{C}(t)$, because $Y_{3}\subset Y$.
Then $\mathrm{rk}\,\mathrm{Pic}(S)\geqslant 2$.

To conclude the~proof we must prove that
$\mathrm{rk}\,\mathrm{Pic}(S)=2$, because there is an exact
sequence
$$
0\longrightarrow\mathbb{Z}\Big[\pi^{*}\big(L\big)\Big]\longrightarrow\mathrm{Pic}\big(\tilde{Y}\big)\longrightarrow \mathrm{Pic}\big(S\big)\longrightarrow 0,%
$$
because  every fibre of $\tau$ is reduced and irreducible (see
the~proof of Proposition~32~in~\cite{BrCoZu03}).

Let $\breve{S}$ be an example of the surface $S$ that is given by
the equation
$$
x\big(q(t)x^{2}+p(t)w^{2}\big)+y\big(r(t)y^{2}+s(t)z^{2}\big)=0\subset\mathrm{Proj}\Big(\mathbb{C}[x,y,z,t]\Big),
$$
where $q(t)$, $p(t)$, $r(t)$, $s(t)$ are~polynomials such that the
inequalities
$$
\mathrm{deg}\big(q(t)\big)>0,\ \mathrm{deg}\big(p(t)\big)\geqslant 0,\ \mathrm{deg}\big(r(t)\big)>0, \mathrm{deg}\big(q(t)\big)\geqslant 0%
$$
hold. The existence of the surface $\breve{S}$ follows from the
equation of the threefold $Y$.

Let $\mathbb{K}$ be an algebraic closure of the field
$\mathbb{C}(t)$, let $L$ be a line $x=y=0$, and let
$$
\gamma\colon\breve{S}\to \mathbb{P}^{1}
$$
be a projection from  $L$. Then $\gamma$ is a conic bundle
defined over $\mathbb{C}(t)$. But $\gamma$ has five~geometrically
reducible fibres $F_{1}, F_{2}, F_{3}, F_{4}, F_{5}$ defined over
$\mathbb{F}$ such that
\begin{itemize}
\item $F_{i}=\tilde{F}_{i}\cup \bar{F}_{i}$, where $\tilde{F}_{i}$ and $\bar{F}_{i}$ are geometrically irreducible curves,%
\item the curve $L\cup F_{i}$ is cut out on the surface
$\breve{S}$ by the equation
$$
y=\epsilon^{i}\sqrt[3]{\frac{q(t)}{r(t)}}x,
$$
where $\epsilon=-(1+\sqrt{-3})/2$ and $i\in\{1,2,3\}$,%
\item the curve $F_{4}\cup L$ is cut out on the surface $\breve{S}$ by the equation $x=0$,%
\item the curve $F_{5}\cup L$ is cut out on the surface $\breve{S}$ by the equation $y=0$.%
\end{itemize}

The~group $\mathrm{Gal}(\mathbb{K}/\mathbb{C}(t))$ naturally acts
on the~set
$$
\Sigma=\Big\{\tilde{F}_{1}, \tilde{F}_{2}, \tilde{F}_{3},
\tilde{F}_{4}, \tilde{F}_{5}, \bar{F}_{1}, \bar{F}_{2},
\bar{F}_{3}, \bar{F}_{4}, \bar{F}_{5}\Big\},
$$
because the conic bundle $\gamma$ is defined over $\mathbb{C}(t)$.
The inequality $\mathrm{rk}\,\mathrm{Pic}(\breve{S})>2$ implies
the~existence of a~subset $\Gamma\subsetneq\Sigma$ consisting of
disjoint curves such that $\Gamma\subsetneq\Sigma$ is
$\mathrm{Gal}(\mathbb{K}/\mathbb{C}(t))$-invariant.

The action of $\mathrm{Gal}(\mathbb{K}/\mathbb{C}(t))$ on the~set
$\Sigma$ is easy to calculate explicitly. Putting
$$
\Delta=\Big\{\tilde{F}_{1}, \tilde{F}_{2}, \tilde{F}_{3}, \bar{F}_{1}, \bar{F}_{2}, \bar{F}_{3}\Big\},\ \Lambda=\Big\{\tilde{F}_{4}, \bar{F}_{4}\Big\},\ \Xi=\Big\{\tilde{F}_{5}, \bar{F}_{5}\Big\},%
$$
we see that the group $\mathrm{Gal}(\mathbb{K}/\mathbb{C}(t))$
acts transitively on each subset $\Lambda$, $\Xi$, $\Delta$,
because we may assume that  $q(t)$, $p(t)$, $r(t)$, $s(t)$ are
sufficiently general. But each subset $\Lambda$, $\Xi$, $\Delta$
does not consist of disjoint curves. Hence, the~equality
$\mathrm{rk}\,\mathrm{Pic}(\breve{S})=2$ holds, which implies that
$\mathrm{rk}\,\mathrm{Pic}(S)=2$.
\end{proof}

The linear system $|\pi^{*}(M-d_{2}L)-E|$ does not have base
points and induces a $\mathbb{P}^{2}$-bundle
$$
\tau\colon
\tilde{V}\longrightarrow\mathrm{Proj}\Big(\mathcal{O}_{\mathbb{P}^{1}}\big(d_{1}\big)\oplus\mathcal{O}_{\mathbb{P}^{1}}\big(d_{2}\big)\Big)\cong\mathbb{F}_{r},
$$
where $r=d_{1}-d_{2}$. Let $l$ be a fibre of the~natural
projection $\mathbb{F}_{r}\to \mathbb{P}^{1}$, and $s_{0}$ be an
irreducible curve on the surface $\mathbb{F}_{r}$ such that
$s_{0}^{2}=r$, and $s_{0}$ is a section of the~projection
$\mathbb{F}_{r}\to \mathbb{P}^{1}$. Then
$$
\pi^{*}\big(M-d_{2}L\big)-E\sim\tau^{*}\big(s_{0}\big)
$$
and $\pi^{*}(L)\sim\tau^{*}(l)$. The morphism $\tau$ induces a
conic bundle
$\tilde{\tau}=\tau\vert_{\tilde{Y}}\colon\tilde{Y}\to\mathbb{F}_{r}$.

Let $\Delta$ be the~degeneration divisor of the conic bundle
$\tilde{\tau}$. Then
$$
\Delta\sim 5s_{\infty}+\mu l,
$$
where $\mu$ is a natural number, and $s_{\infty}$ is
the~exceptional section of the surface $\mathbb{F}_{r}$.

Let $S$ be a surface in $\tilde{Y}$ and $B$ be a threefold in
$\tilde{V}$ dominating the~curve $s_{0}$. Then
$$
B\cong\mathrm{Proj}\Big(\mathcal{O}_{\mathbb{P}^{1}}\big(d_{1}\big)\oplus\mathcal{O}_{\mathbb{P}^{1}}\big(d_{3}\big)\oplus\mathcal{O}_{\mathbb{P}^{1}}\Big)
$$
and $\pi(B)\cong B$. But $\pi(B)\cap Y=\pi(S)\cup Y_{3}$.

The surface $Y_{3}$ is cut out on $\pi(B)$ by the equation
$x_{1}=0$, where $\pi(B)\in |M-d_{2}L|$. We have
$$
S\sim 2T+\big(d_{1}+n\big)F,
$$
where $T$ is a~tautological line bundle on $B$, and $F$ is a~fibre
of the~projection $B\to \mathbb{P}^{1}$. Then
$$
K_{S}^{2}=-5d_{1}+2d_{3}-4d_{2}-3n+8
$$
and
$\mu=s_{0}\cdot\Delta=5d_{1}-2d_{3}+4d_{2}+3n$.

It follows from the equivalence $2K_{\mathbb{F}_{r}}+\Delta\sim
s_{\infty}+(3d_{1}-2d_{3}+6d_{2}+3n-4)l$ that
$$
\big|2K_{\mathbb{F}_{r}}+\Delta\big|\ne\emptyset\iff 3d_{1}-2d_{3}+6d_{2}+3n\geqslant 4,%
$$
which implies that $Y$ is nonrational by Theorem~10.2 in
\cite{Sho83} if $3d_{1}-2d_{3}+6d_{2}+3n\geqslant 4$.

The threefold $Y$ is nonruled if and only if it is nonrational,
because the threefold $Y$ is rationally connected. So,
the~threefold $X$ is nonrational by Theorem~1.8.3 in \S IV of
the~book \cite{Ko96} whenever
$$
3d_{1}-2d_{3}+6d_{2}+3n\geqslant 4,
$$
which implies that $3d_{1}-2d_{3}+6d_{2}+3n<4$, because we assume
that $X$ is rational.

We see that either $d_{1}=0$ and $n=1$ or $d_{1}=1$ and
$d_{2}=n=0$ by Lemmas~\ref{lemma:smoothness},
\ref{lemma:smoothness-extra} and \ref{lemma:smoothness-rkPic-2},
but the~threefold $X$ is birational  to a smooth cubic threefold
in the case when $d_{1}=1$~and~$d_{2}=n=0$, which is nonrational
by \cite{ClGr72}. Then $d_{1}=0$ and $n=1$. The~assertion of
Theorem~\ref{theorem:main} is~proved.

\end{document}